\documentclass[11pt,twoside]{article}
\usepackage{amssymb}

\def\flex#1{\mathrel{\mathop{\kern 0pt\hbox to 
10mm{\rightarrowfill}}\limits_{#1
\rightarrow \infty}}}

\def\Flex#1{\mathrel{\mathop{\kern 0pt\hbox to 
10mm{\rightarrowfill}}\limits_{#1
\rightarrow \infty}}}

\font\bb=msbm10 at 12pt

\def\Q{\hbox{\bb Q}}
\def\R{\hbox{\bb R}}

\def\P{\hbox{\bb P}}

\def\I{{\bf I}}
\def\J{{\bf J}}

\def\1{1\hspace{-1.2mm}\mbox{{\normalsize I}}}

\def\Dem{\vskip 4mm\noindent { Proof{}.}  }

\newtheorem{Th}{Theorem}

\newtheorem{Lemme}[Th]{Lemma}

\newtheorem{Prop}[Th]{Proposition}

\title{Boundary behavior of a constrained Brownian motion between reflecting-repellent walls}
\author{Dominique L\'epingle}

\begin{document}
\maketitle

\begin{abstract}
Stochastic variational inequalities provide a unified treatment for stochastic differential 
equations living in a closed domain with normal reflection
 and (or) singular repellent drift. When the domain is a polyhedron, we prove that the reflected-repelled 
Brownian motion does not hit the non-smooth part of the boundary. A sufficient condition for non-hitting a 
face of the polyhedron is derived from the one-dimensional case. A complete answer to the question of
attainability of the walls of the Weyl chamber may be given for a radial Dunkl process.
\end{abstract} 
\section{Introduction}

There have been many works about stochastic differential equations with
reflection on the boundary of a domain. In some of them the domain is a convex polyhedron (\cite{HR81}, \cite{VW85}, \cite{W87}, \cite{DW96}, \cite{De08}]. 
A typical question in this setting is the following: does the continuous
process hit the non-smooth
part of the boundary? The answer depends on the drift and diffusion coefficients of the process and on the direction of reflection (normal or 
oblique). In particular, R.Williams \cite{W87} has proven that the 
Brownian motion with a skew symmetry condition on 
the direction of reflection does not touch the intersections of the faces
of the polyhedron. 

On the other hand there also exists an extensive literature about 
non-colliding Brownian particles (\cite{D62}, \cite{B91}, \cite{HW96}, 
\cite{G99}, \cite{O03}). Most of these 
works originate in the study of the eigenvalues of Gaussian matrix processes.
These eigenvalues are solutions to systems of stochastic differential equations with a singular drift that prevents the particles from colliding. Extensions of 
these systems are Dunkl processes \cite{RV98} that have recently been developed 
in connection with harmonic analysis on symmetric spaces. The radial part
of a Dunkl process may be considered as a Brownian motion perturbed by a singular drift which forces the process to live in a cone generated by the 
intersection of a finite set of half-spaces (\cite{C06}, \cite{CGY08}).
Depending on the values of 
some parameter, the process may touch the walls of the cone or not.

Actually it is possible to unify both theories of (normal) reflection and 
strong repulsion
within a common framework. This is the role of stochastic variational inequalities, also called multivalued stochastic differential equations 
(MSDE) that were mainly developed by E.C\'epa (\cite{C95}, \cite{C98}). 
These equations are associated to a convex function 
in a domain of $\R^d$. Depending on the boundary behavior of this
function the diffusion will (normally) reflect on the boundary, hit the boundary without local time, or
live in the open domain. We shall here follow this way and 
concentrate on a Brownian motion living in a convex polyhedral domain, bounded
 or unbounded. To each face of
 the polyhedron is associated a repelling force with normal reflection when the repulsion is not
 strong enough. In this setting we shall ask whether the process may 
hit the various faces. Our first
 task will be to rule out the possibility of hitting the intersection of two faces. Once this
 is achieved, the problem is now basically one-dimensional and we may use the ordinary 
 scale function of real diffusions.
 
 In several previous works (\cite{I06}, \cite{CL07}), this issue has been studied in the particular case of
 the hyperplanes $H_{ij}:=\{x=(x_1,\ldots,x_d)\in\R^d:\,x_i=x_j\}, \,i\neq j$ 
and presented as the 
problem of collisions between Brownian particles. There is a simple collision if two coordinates
coincide and a multiple collision if at least three coordinates coincide at the same time. 
Because
 the $d$-dimensional Brownian motion does not hit the intersection of two hyperplanes, one can
 guess that an additional drift does not change anything. However a rigorous proof is necessary 
because the singularity of the drift makes useless the usual Girsanov change of probability
measure. The counterexample of Bass and Pardoux \cite{BP87} also showed that uniform nondegeneracy 
of the diffusion term does not preclude multiple collisions.

As in \cite{CL07} where the particular case of 
electrostatic repulsion was considered, our proof only 
uses basic tools from stochastic calculus, mainly McKean's 
martingale method \cite{M69} which was 
already used in \cite{B89} to prove non-collision for the eigenvalues of Wishart
 processes. Another way could be to use the theory of Dirichlet forms as done in 
\cite{I06} where a general condition of non-collision has been obtained.

The paper is organized as follows. In Section 2 we introduce basic definitions and notations.
The main features about stochastic variational inequalities are also recalled. Section 3 is
devoted to non attainability of the edges of the polyhedron. 
In Section 4 we give a sufficient condition of non attainability of a single face. Section 5 presents some
 applications to Brownian particles with nearest neighbor interaction, Wishart processes and 
Dunkl processes. 

\section{Multivalued stochastic differential equation in a polyhedral domain}

Let $(\Omega,{\cal F},({\cal F}_t,t\geq 0),\P)$ be a filtered probability space endowed
with the usual conditions
and $B=(B_t)$ be a $({\cal F}_t)$-adapted $d$-dimensional Brownian motion starting from the
origin. Let 
\begin{equation} 
  \Phi : \R^d \rightarrow (-\infty,+ \infty]
\end{equation}
be a lower semi-continuous convex function such that 
\begin{equation}
  dom(\Phi):=\{x\,:\, \Phi(x)\,<\,\infty\}
\end{equation}
has nonempty interior. Let 
\begin{equation}
  D:=Int(dom(\Phi))\;.
\end{equation}
For simplicity of notation, we will assume that $\Phi$ is $C^1$ on $D$. If $x\in \partial D$,
we say that the unit vector $n(x)$ is a unit inward normal to $D$ at $x$ if 
\begin{equation}
  n(x). (x-z) \,\leq 0
\end{equation}
for any $z\in \overline{D}$. Based on the results in \cite{C95}, the following theorem has been proved 
in \cite{CL97} (see also Theorem 2.2 in \cite{CL01}).
\begin{Th} 
For any  ${\cal F}_0$-measurable random variable $X_0$ with values in $\overline{D}$, there exist a unique
continuous $({\cal F}_t)$-adapted process $X=\{X_t,0\leq t < \infty\}$ with values in 
$\overline{D}$ and a unique continuous $({\cal F}_t)$-adapted non-decreasing process 
$L=\{L_t,0\leq t < \infty\}$ such that
\begin{equation}
 \begin{array}{llll}
   X_t & = & X_0 \, +\,B_t\,- \, \int_0^t \nabla \Phi(X_s)\,ds\,+\,\int_0^t n_s\,dL_s
        &\qquad 0\leq t<\infty \\
   L_t & = & \int_0^t{\bf 1}_{\{X_s\in\partial D\}}dL_s &\qquad 0\leq t<\infty 
 \end{array}
\end{equation}
where $n_s$ is $dL_s$-a.e. a unit inward normal to $D$ at $X_s$. For any $0<T<\infty$, 
\begin{equation}
  \int_0^T {\bf 1}_{\{X_s\in\partial D\}} ds \,=\,0
\end{equation}
and
\begin{equation}
 \int_0^T |\nabla \Phi(X_s)|\,ds\,<\,\infty \; .
\end{equation}
\end{Th}

From now on we concentrate on a particular polyhedral setting. Let ${\bf I}:=\{1,\ldots,m\}$ where
$m\geq 1$.
We consider a convex function $\Phi$ of the following form
\begin{equation}
 \Phi(x)\,:=\, \sum_{i\in {\bf I}}\phi_i(x.n_i-a_i)
\end{equation}
where for any $i\in {\bf I}$,
\begin{equation}
 \begin{array}{ll}
 \phi_i & \mbox{ is a convex l.s.c. function },\; \;\phi_i=+\infty \; \mbox{ on }(-\infty,0), \;\;
   \phi_i \mbox{ is } C^1 \mbox{ on }(0,+\infty) \\
 n_i & \mbox{ is a unit vector} \\
 a_i & \mbox{ is a real number}\;. 
 \end{array}
\end{equation} 
We may assume all $n_i$ are different. Then,
\begin{equation}
\begin{array}{lll}
  \nabla\Phi(x) & = & \sum_{i\in {\bf I}} n_i\,\phi_i^{\prime} (x.n_i-a_i) \\
  D & = & \{x\in \R^d: \,x.n_i\,>\,a_i  \;\;\forall i\in {\bf I} \} \\
  \overline{D} & = & \{x\in \R^d: \,x.n_i\,\geq \,a_i \;\;\forall i\in {\bf I} \} \;.  
\end{array}
\end{equation}
As $D$ is not empty, there exists a ball with center $y\in D$ and radius $b>0$
 included in $D$. Let $X_t$ be the solution given by Theorem 1. For $i\in {\bf I}$ let
\begin{equation}
 U_t^i :=X_t.n_i\,-\,a_i \;.
\end{equation}
We will need a strengthening of inequality (7) (\cite{CL01},Th.2.2). 
\begin{Lemme}
For any $i\in I$, for any $0<t<\infty$,
\begin{equation}
 \int_0^t|\phi_i^{\prime}(U_s^i)|\,ds \,<\,\infty\;.
\end{equation} 
\end{Lemme}
\Dem
This is clear if $\phi_i^{\prime}(0+)>-\infty$. Let
\begin{equation}
  {\bf J}:=\{j\in {\bf I}:\phi_j^{\prime}(0+)=-\infty \}
\end{equation}
and let $0<\varepsilon<b$ be such that $\phi_j^{\prime}(u)<0$ 
for any $j\in {\bf J}$ and $u\in (0,\varepsilon)$.
For ${\bf K}\subset {\bf J}$ let
\begin{equation}
 A_{\bf K}:= \{ x\in \overline{D}: x.n_j<a_j+\varepsilon \;\;\forall j\in {\bf K},\;
  \;\; x.n_j \geq a_j+\varepsilon \;\; \forall j\in {\bf J}\setminus {\bf K} \} \;.
\end{equation}
Then for $t>0$
\begin{equation}
\int_0^t{\bf 1}_{A_{\bf K}}(X_s)|\sum_{j\not\in {\bf K}}n_j\,\phi_j^{\prime}(U_s^j)|ds
  \leq \sum_{j\not\in {\bf K}}\int_0^t{\bf 1}_{A_{\bf K}}(X_s) |\phi_j^{\prime}(U_s^j)|ds\,<
\,\infty \;.
\end{equation}
Using (7) we get
 \begin{equation}
\int_0^t{\bf 1}_{A_{\bf K}}(X_s)|\sum_{j\in {\bf K}}n_j\,\phi_j^{\prime}(U_s^j)| ds \,<\,\infty
\end{equation}
and therefore
\begin{equation}
 \begin{array}{lll}
  -(b-\varepsilon)\sum_{j\in {\bf K}}\int_0^t{\bf 1}_{A_{\bf K}}(X_s)\phi_j^{\prime}(U_s^j)\, ds
    & \leq  & \int_0^t{\bf 1}_{A_{\bf K}}(X_s)
\sum_{j\in {\bf K}}(y-X_s).n_j|\phi_j^{\prime}(U_s^j)| ds \\
    & \leq  & \int_0^t{\bf 1}_{A_{\bf K}}(X_s)|y-X_s|
\sum_{j\in {\bf K}}|\phi_j^{\prime}(U_s^j)| ds \\
 & < & \infty 
  \end{array}
 \end{equation}
from the continuity of $X$ on $[0,t]$. Then for any $j\in {\bf J}$ 
\begin{equation}
 \begin{array}{lll}
 \int_0^t|\phi_j^{\prime}(U_s^j)|\,ds  & = & \int_0^t{\bf 1}_{\{U_s^j<\varepsilon\}}
  |\phi_j^{\prime}(U_s^j)|\,ds \,+\,  \int_0^t{\bf 1}_{\{U_s^j\geq \varepsilon\}}
  |\phi_j^{\prime}(U_s^j)|\,ds \\
   & = &\sum_{j\in {\bf K}\subset {\bf J}}\int_0^t{\bf 1}_{A_{\bf K}}(X_s)\,
  |\phi_j^{\prime}(U_s^j)| ds \,+\,  \int_0^t{\bf 1}_{\{U_s^j\geq \varepsilon\}}
  |\phi_j^{\prime}(U_s^j)|\,ds \\
   & < & \infty \; . \hfill \square
  \end{array}   
\end{equation} 
 
For any ${\bf J}\subset {\bf I}$, ${\bf J}\neq \emptyset$, we set
\begin{equation}
  \begin{array}{lll}
   H_{{\bf J}} & := & \{ x\in \R^d : x.n_j=a_j \;\, \forall j\in {\bf J}\} \\
   K_{{\bf J}} & := & \{ x\in \R^d : x.n_j=a_j \;\, \forall j\in {\bf J},\;\; x.n_j>a_j \;\,\forall j\not\in {\bf J}\} \\
   \sigma_{{\bf J}} & := & \inf \{t>0 \,:\, X_t\,\in H_{{\bf J}}\} \\ 
   \tau_{{\bf J}} & := & \inf \{t>0 \,:\, X_t\,\in K_{{\bf J}}\} \ \;.
  \end{array}
\end{equation}
 
\begin{Lemme}
Let ${\bf J}\subset {\bf I}$ and $V:=span\{n_j,j\in {\bf J}\}$. If $n(x)$ is a 
unit inward normal to $D$ at $x\in K_{{\bf J}}$, then $n(x)\in V$.
\end{Lemme}
\Dem
Let $v\perp V$. For  $\varepsilon >0$ small enough, 
\[
  \begin{array}{lll}
    z_1=x+ \varepsilon v & &
    z_2=x- \varepsilon v
  \end{array}
\]
satisfy 
\[
  \begin{array}{lll}
   z_1.n_j=a_j \;\; \forall j\in {\bf J}  & & \;\; z_1.n_i>a_i \;\; \forall i\not\in {\bf J}  \\
   z_2.n_j=a_j \;\; \forall j\in {\bf J}  &  & \;\; z_2.n_i>a_i \;\; \forall i\not\in {\bf J} \;. 
  \end{array}
\]
Then 
\[
  \begin{array}{lll}
   n(x).(x-z_1) \leq 0 & & n(x).(x-z_2) \leq 0 
  \end{array}
\]
and therefore
\[
  n(x).v=0  \;. 
\] 
$\hfill \square$

\section{Nonattainability of the edges}

This section is devoted to the proof of the following theorem.
\begin{Th}
For any ${\bf J}\subset {\bf I}$ with $card({\bf J})\geq 2$,
\[
\P(\sigma_{{\bf J}}  = \infty)\,=\, 1 \;.
\]
\end{Th}
\Dem
a/ We first consider the initial condition $X_0$. From (6) we deduce that  
for any $u>0$ there exists $0<v<u$ such that $X_v\in D$ a.s. Using the 
continuity of paths and the Markov property we may and do assume that 
$X_0\in D$ in order to prove that $\sigma_{{\bf J}}=\infty$ a.s.

b/ We will also assume that 
\begin{equation}
  \max_{i\in {\bf I}} \phi_i^{\prime}(0+) \,<\,0 \;.
\end{equation}
If not we introduce for any $0<T<\infty$ the equivalent probability measure $\Q$ defined on ${\cal F}_T$ by 
\[
  \frac{d\Q}{d\P} := \exp \{c(B_T.\sum_{i\in {\bf I}}n_i)-\frac{1}{2}c^2 T
     |\sum_{i\in {\bf I}}n_i|^2 \}
\]
where 
\[
  c> \max_{i\in {\bf I}} \phi_i^{\prime}(0+)\;.
\]
The continuous process 
\[
  B_t^{\prime}:=B_t-ct\sum_{i\in {\bf I}}n_i
\]
is a $\Q$-Brownian motion on $[0,T]$ and now
\[
  dX_t=dB_t^{\prime}-\sum_{i\in {\bf I}}n_i
  \psi_i^{\prime}(X_t.n_i-a_i)dt-n_tdL_t
\]
where 
\[
  \psi_i(u):=\phi_i(u)-cu \quad i\in{\bf I}\;.
\]  
If $\Q(\sigma_{{\bf J}}<T)=0$ then $\P(\sigma_{{\bf J}}<T)=0$ 
and if this is true
for any $T$ we obtain $\P(\sigma_{{\bf J}}=\infty)=1$.

c/ We are now going to prove that $\sigma_{{\bf I}}=\tau_{{\bf I}}=\infty$ a.s.
(with $m\geq 2$). For any ${\bf J}\subset {\bf I}$ let
\begin{equation}
 \begin{array}{lll}
   V_{\J} & := & span\{n_j,j\in \J\}\\
   q_{\J} & := & dim\, V_{\J}\\
   \pi_{\J} & := & \mbox{orthogonal projection onto }\;V_{\J}\;.
 \end{array}
\end{equation} 
If $q_{\I}=1$, then $m=2$, $n_1+n_2=0$ and $H_{\I}=K_{\I}=\emptyset$. Assume now
$ q_{\I}\geq 2$ and $H_{\I}\neq \emptyset$. Choose some $z\in H_{\I}$ and set
\begin{equation}
  Z_t:= \pi_{\I}(X_t-z) \;.
\end{equation}  
Then
\begin{equation}
  Z_t\,=\,Z_0\,+\,C_t \,-\sum_{i\in \I}\int_0^tn_i\,\phi_i^{\prime}(U_s^i)ds\,
     +\, \int_0^tn_s\,dL_s
\end{equation}
where $C$ is a $q_{\I}$-dimensional Brownian motion. Set
\[
  S_t\,:=\,|Z_t|^2\;.
\]
Then 
\begin{equation}
 S_t=S_0 + 2\int_0^t Z_s.dC_s -2\sum_{i\in \I}\int_0^t U_s^i\,\phi_i^{\prime}
  (U_s^i)ds +2\int_0^t Z_s.n_s\,dL_s + q_{\I}t .
\end{equation}  
From Lemma 2 we deduce that on $\partial D=\cup_{\J\subset\I}K_{\J}$
\[
  Z_s.n_s\,=\,(X_s-z).n_s\,=\,0
\]
and thus
\[
  \int_0^t Z_s.n_s\,dL_s = 0\;.
\]  
Let $0<T<\infty$. For $t<\tau_{\I}\wedge T$,
\begin{equation}
 \log S_t \,=\, \log S_0 \,+\, 2\int_0^t \frac{Z_s.dC_s}{S_s} -2 \sum_{i\in \I}\int_0^t \frac{U_s^i \, \phi_i^{\prime}
  (U_s^i)}{S_s}ds \,+\,(q_{\I}-2)\int_0^t \frac{ds}{S_s} \;.
\end{equation}
  From the assumption made in b/ there exists $0<c\leq \infty$ such that $\phi_i^{\prime}\leq 0$ on $(0,c]$ and
\begin{equation}
  \begin{array}{lll}
   -\int_0^t \frac{U_s^i \,\phi_i^{\prime}(U_s^i)}{S_s}ds  & \geq &
    -\int_0^t \frac{U_s^i \,\phi_i^{\prime}(U_s^i)}{S_s}\,{\bf 1}_{\{U_s^i\geq c\}}ds \\
       & \geq & -\frac{1}{c} \int_0^T |\phi_i^{\prime}(U_s^i)|\,ds \\
   & > & -\infty  \;.
  \end{array}
\end{equation}
We now proceed as in (\cite{M69},p.47).
As $t\rightarrow \tau_{\I}\wedge T$, the local martingale part in the r.h.s. of
(25) either converges 
to a finite limit or oscillates between $+\infty$ and $-\infty$. Thus it does not 
converge to $-\infty$ and a.s. $S_{\tau_{\I}\wedge T}>0$. Therefore
\[
  \P(\tau_{\I} \leq T)=0
\]
and the conclusion follows since $T$ is arbitrary.

d/ Let now $\J \subset \I $ with $2\leq |\J| \leq m-1$.
We shall show by a backward induction on $|\J|$ that $\P(\tau_{\J}=\infty)=1$. Remark that the 
backward induction assumption entails the equality $\sigma_{\J}= \tau_{\J}$ a.s.. 
As previously done we may assume $q_{\J} \geq 2$ and $K_{\J}\neq \emptyset$. Select now 
$z\in K_{\J}$ and set
\begin{equation}
\begin{array}{lll}
  Z_t &:= & \pi_{\J}(X_t-z) \\
   & = & Z_0\,+\,C_t \,-\sum_{j\in \J}\int_0^tn_j\,\phi_j^{\prime}(U_s^j)ds
     -\sum_{i\not \in \J}\int_0^t \pi_{\J}n_i\,\phi_i^{\prime}(U_s^i)ds
     +\, \int_0^t \pi_{\J}n_s\,dL_s 
 \end{array}
\end{equation} 
where $C$ is a $q_{\J}$-dimensional Brownian motion. Let again $S_t:=|Z_t|^2$. For 
$\varepsilon>0$ and $r>0$ we set
\begin{equation}
 \begin{array}{lll}
   \tau_{\varepsilon} & := & \inf\{t>0\,:\,S_t+\min_{i\not\in \J}(U_t^i)^2
    \leq \,2\,\varepsilon^2 \} \\
   \rho_{r}  & = & \inf\{t>0 \,:\,|X_t|\geq r\} \;.
 \end{array}
\end{equation}
 From the induction assumption we infer that $\tau_{\varepsilon}\rightarrow \infty$ as
$\varepsilon$ goes to $0$. Let $0<T<\infty$. We introduce the equivalent probability measure 
$\Q$ defined on ${\cal F}_T$ by 
\begin{equation}
\begin{array}{lll}
  \frac{d\Q}{d\P}
& = & \exp \{ \int_0^ {\tau_{\varepsilon} \wedge \rho_{r}\wedge T}
   \sum_{i\not\in \J}{\bf 1}_{\{U_s^i\geq \varepsilon\}}\phi_i^{\prime}(U_s^i)\,n_i.dC_s \\
  &  & \qquad -\frac{1}{2}\int_0^{\tau_{\varepsilon} \wedge \rho_{r}\wedge T}
    | \sum_{i\not\in \J}{\bf 1}_{\{U_s^i\geq \varepsilon\}}\phi_i^{\prime}(U_s^i)\,
     \pi_{\J}n_i|^2\,ds\} \; .
 \end{array}
 \end{equation}  
Then
\[
  D_t \,:=\,C_t- \int_0^ {\tau_{\varepsilon} \wedge \rho_{r} \wedge T}
   \sum_{i\not\in \J}{\bf 1}_{\{U_s^i\geq \varepsilon\}}\phi_i^{\prime}(U_s^i)\pi_{\J}n_i\,ds
\]
is a $q_{\J}$-dimensional $\Q$-Brownian motion on $[0,T]$. For $t\leq 
\tau_{\varepsilon} \wedge \rho_{r}\wedge T$,
\begin{equation}
 \begin{array}{lll}
  S_t & = & S_0+2\int_0^t Z_s.dD_s -2\sum_{i\in \I}\int_0^t U_s^i\,\phi_i^{\prime}
  (U_s^i)ds \,- \,2\sum_{i\not\in\J}
\int_0^t {\bf 1}_{\{U_s^i<\varepsilon\}}\,Z_s.n_i\,
  \phi_i^{\prime}(U_s^i)\,ds \\
   &  & +2\sum_{{\bf L}\subset\I,{\bf L}\not\subset\J}\int_0^t{\bf 1}_{K_{{\bf L}}}(X_s)\,Z_s.n_s dL_s 
       +q_{\J}t
 \end{array}
\end{equation}
 and for $t<\sigma_{\J}\wedge \tau_{\varepsilon}\wedge \rho_{r} \wedge T$,
\begin{equation}
 \begin{array}{lll}
 \log S_t & = & \log S_0 + 2\int_0^t \frac{Z_s.dD_s}{S_s} -2 \sum_{j\in \J}\int_0^t 
 \frac{U_s^j \, \phi_j^{\prime}(U_s^j)}{S_s}ds  \\
   &  &- 2 \sum_{i\not\in \J}\int_0^t {\bf 1}_{\{U_s^i<\varepsilon\}}
 \frac{\phi_i^{\prime}(U_s^i)}{S_s} Z_s.n_i\,ds \\
& & + 2\sum_{{\bf L}\subset\I, {\bf L}\not\subset\J}\int_0^t {\bf 1}_{K_{{\bf L}}}(X_s)
  \frac{Z_s.n_s}{S_s}dL_s  \\
   & & +(q_{\J}-2)\int_0^t \frac{ds}{S_s} \;.
 \end{array}
\end{equation}     
From the induction hypothesis and the continuity of paths, if $\sigma_{\J}
<\infty$ for any ${\bf L}\not\subset \J$ there exists an interval $(\sigma_{\J}-\delta,\sigma_{\J}]$
of positive length on which $X_s\not\in K_{{\bf L}}$. Therefore
\begin{equation}
-\int_0^{\sigma_{\J}\wedge \tau_{\varepsilon} \wedge \rho_{r} \wedge T}
{\bf 1}_{K_{{\bf L}}}(X_s)\frac{Z_s.n_s}{S_s}dL_s >-\infty \;.
\end{equation}
For $s<\tau_{\varepsilon}$, if $U_s^i<\varepsilon$ for some $i\not\in\J$, then
$S_s\geq \varepsilon^2$ and we obtain as well
\begin{equation}
  -\int_0^{\sigma_{\J}\wedge \tau_{\varepsilon}\wedge \rho_{r} \wedge T}
  {\bf 1}_{\{U_s^i<\varepsilon\}}
 \frac{\phi_i^{\prime}(U_s^i)}{S_s} Z_s.n_i\,ds >-\infty\;.
\end{equation}
The other terms behave as in c/ and thus
\begin{equation}
  0= \Q(\sigma_{\J}\leq \tau_{\varepsilon} \wedge \rho_{r} \wedge T)
    = \P(\sigma_{\J}\leq \tau_{\varepsilon} \wedge \rho_{r} \wedge T) \;.
\end{equation}
Letting $\varepsilon$ go to $0$, $r$ and $T$ to $\infty$ we get
\[
   \P(\sigma_{\J}=\infty)=1 
\]
and we are done. $\hfill \square$   

\section{Keeping off from a wall}

We first recall some facts in the one-dimensional setting \cite{LM87}. Let 
$\phi:\R\rightarrow (-\infty,+\infty]$ be
a convex lower semicontinuous function.
Assume $\phi=+\infty$ on $(-\infty,0)$ and $C^1$ on $(0,+\infty)$. 
Consider the one-dimensional equation
\begin{equation}
 \begin{array}{lll}
   dY_t &= & dB_t-\phi^{\prime}(Y_t)dt+\frac{1}{2}dL^0_t \\
   Y_t & \geq & 0
 \end{array}
\end{equation}  
where $L^0$ is the local time of $Y$ at $0$. There are three types of boundary
behavior: 

\vspace{0.3cm}

 \begin{tabular}{|l|l|}  \hline
   &  repulsion \\  \hline
 $\phi(0)<\infty$  & weak: local time not zero \\ \hline
 $\phi(0)=\infty , \,\int_{0+}\exp\{2\phi\}<\infty $
   &  middle: local time zero \\ \hline
 $\phi(0)=\infty ,\,\int_{0+}\exp\{2\phi\}=\infty $
    &  strong: boundary not hit \\ \hline 
 \end{tabular}
 
\vspace{0.3cm }

We shall check the behavior of the multidimensional process $X$ accords with this
 classification in the neighborhood of the faces of the polyhedron.
For any $i\in \I$ we respectively write $H_i,K_i,\sigma_i,\tau_i$  
in place of $H_{\{i\}},K_{\{i\}},\sigma_{\{i\}},\tau_{\{i\}}$.

\begin{Prop}
For any $i\in\I$ such that $\phi_i(0)=\infty$ and any $t>0$,
\begin{equation}
  \int_0^t{\bf 1}_{H_i}(X_s)\,dL_s\,=\,0\;.
\end{equation}
\end{Prop}
\Dem  
From the occupation times formula and Lemma 1 we obtain
\begin{equation}
  \int_0^{\infty}
  L_t^a(U_i)\,|\phi_i^{\prime}(a)|\,da = \int_0^t|\phi_i^{\prime}(U_s^i)|\,ds
    <\infty
    \end{equation}
and from $\phi_i(0)=\infty$ and the continuity of 
$a \mapsto L_t^a(U_i)$ we deduce
\begin{equation}
L_t^0(U_i)=0\;.
\end{equation}
Thus
\begin{equation}
\begin{array}{lll}
 0 &= &U_t^i-(U_t^i)^+ \\
  &= &\int_0^t {\bf 1}_{H_i}(X_s)n_i.dB_s - \int_0^t{\bf 1}_{H_i}(X_s)
  \sum_{j\in \I}\phi_j^{\prime}(U_s^j)\,n_i.n_j\,ds +\int_0^t
  {\bf 1}_{H_i}(X_s)\,n_i.n_s\,dL_s \\
   & = & \int_0^t{\bf 1}_{K_i}(X_s)\,n_i.n_s\,dL_s  \\
   & = & \int_0^t{\bf 1}_{K_i}(X_s)\,dL_s  \\
   & = & \int_0^t{\bf 1}_{H_i}(X_s)\,dL_s \; .\hfill \square
\end{array}
\end{equation}
 
We now set for any $i\in\I$ and $x\geq 0$
\[
    p_i(x):=\int_1^x\exp\{2(\phi_i(u)-\phi_i(1))\}\,du \;.
\]

\begin{Th}
For any $i\in\I$ such that $p_i(0)=-\infty$ or equivalently
\begin{equation}
   \int_{0+}\exp\{2\phi_i\}=\infty \;,
\end{equation}
then $\P(\sigma_i=\infty)=\P(\tau_i= \infty)=1$.
\end{Th}
\Dem  
 From Ito formula and Proposition 5 we obtain 
 \begin{equation}
  p_i(U_t^i)=p_i(U_0^i)+\int_0^t p_i^{\prime}(U_s^i)[dC_s^i-\sum_{j\neq i}
n_i.n_j\,\phi_j^{\prime}(U_s^j)ds
 +\sum_{j\neq i}{\bf 1}_{K_j}(X_s)\,n_i.n_j\,dL_s] 
 \end{equation}
 where $C^i=B.n_i$ is a one-dimensional Brownian motion.
 As in the proof of Theorem 2, let
 \begin{equation}
  \begin{array}{lll}
  \tau_{\varepsilon} & := & \inf\{t>0\,:\,U_t^i+\min_{j\neq i}(U_t^j)
    \leq \,2\,\varepsilon \} \\
   \rho_{r}  & = & \inf\{t>0 \,:\,|X_t|\geq r\} \;.
  \end{array}
\end{equation}
 Let $0<T<\infty$. We again introduce the equivalent probability measure 
$\Q$ defined on ${\cal F}_T$ by 
\begin{equation}
\begin{array}{lll}
  \frac{d\Q}{d\P}
& = & \exp \{ \int_0^ {\tau_{\varepsilon} \wedge \rho_{r}\wedge T}
   \sum_{j\neq i}{\bf 1}_{\{U_s^j\geq \varepsilon\}}\phi_j^{\prime}(U_s^j)\,n_i.n_j \, dC_s^i \\
  &  & \qquad -\frac{1}{2}\int_0^{ \tau_{\varepsilon} \wedge \rho_{r}\wedge T}
    | \sum_{j\neq i}{\bf 1}_{\{U_s^j\geq \varepsilon\}}\phi_j^{\prime}(U_s^j)\,
     n_i.n_j|^2\,ds\} \; .
 \end{array}
 \end{equation}
Then
\begin{equation}
 D_t^i := C_t^i-  
 \int_0^ {t\wedge \tau_{\varepsilon} \wedge \rho_{r}}
\sum_{j\neq i}{\bf 1}_{\{U_s^j\geq \varepsilon\}}\phi_j^{\prime}(U_s^j)\,n_i.n_j\,ds 
\end{equation}
is a $\Q$-Brownian motion on $[0,T]$ and for 
 $t\leq  \tau_{\varepsilon} \wedge \rho_{r} \wedge T$,
\begin{equation}
  p_i(U_t^i)=p_i(U_0^i)+\int_0^t p_i^{\prime}(U_s^i)[dD_s^i-\sum_{j\neq i}
   {\bf 1}_{\{U_s^i<\varepsilon\}} n_i.n_j\,\phi_j^{\prime}(U_s^j)ds
 +\sum_{j\neq i}{\bf 1}_{K_j}(X_s)\,n_i.n_j\,dL_s] \;.
 \end{equation}  
 As in the proof of Theorem 2, for any $j\neq i$,
\begin{equation}
  - \int_0^{\sigma_i\wedge \tau_{\varepsilon} \wedge \rho_{r}\wedge T}
   {\bf 1}_{\{U_s^i<\varepsilon\}}\,p_i^{\prime}(U_s^i)\, n_i.n_j\,\phi_j^{\prime}(U_s^j) \,ds\,> - \infty
\end{equation}
and
\begin{equation}
  +\int_0^{\sigma_i\wedge \tau_{\varepsilon} \wedge \rho_{r}\wedge T}
  {\bf 1}_{K_j}(X_s)\,  p_i^{\prime}(U_s^i) \,n_i.n_j\,dL_s \,>-\infty
\end{equation}
and then
\begin{equation}
 0\,=\,\Q(\sigma_i\leq  \tau_{\varepsilon} \wedge \rho_{r}\wedge T)
 = \,\P(\sigma_i\leq  \tau_{\varepsilon} \wedge \rho_{r}\wedge T)
\end{equation}
meaning that $\P(\sigma_i=\infty)=1$. \hfill $\square$

\section{Applications}

\subsection{Brownian particles with nearest neighbor repulsion}

H.Rost and M.E.Vares \cite{RV85} have considered the following system:
\begin{equation}
 \begin{array}{llll}
  dX_t^1 & = & dB_t^1 \,+\, \phi^{\prime}(X_t^2-X_t^1)\,dt & \\
  dX_t^i & = & dB_t^i\,+\,(\phi^{\prime}(X_t^{i+1}-X_t^i)
     -\phi^{\prime}(X_t^i-X_t^{i-1}))\,dt & \quad i=2,\ldots,n-1 \\
  dX_t^n & = & dB_t^n \,-\, \phi^{\prime}(X_t^n-X_t^{n-1})\,dt 
 \end{array}
\end{equation}
where $X_t^1<\ldots<X_t^n$ and $\phi$ is a positive convex function 
on $(0,\infty)$ satisfying
\begin{equation}
  \phi(0)=\infty \, ,\qquad\phi(\infty)=0 \,,\qquad \int_0^1
   (\phi^{\prime}(x))^2e^{-2\phi(x)}\,dx\,<\,\infty \;. 
\end{equation}
This is a MSDE where function $\Phi$ is given by (8)
with $\phi_i(x)=\phi(\sqrt{2}\,x)$, $n_i=\frac{1}{\sqrt{2}}(e_{i+1}-e_i)$,
$a_i=0$ for $i=1,\ldots,n-1$ and 
$e_j$ the $j$-th basis vector.  
Condition (50) for non-collision is stronger than (40) as can be seen from
Schwarz inequality:
\[
  \infty \,=\,(\phi(0)-\phi(1))^2\,\leq \, \int_0^1(\phi^{\prime})^2
  \, e^{-2\phi}\,\int_0^1\,e^{2\phi} \,.
\]  

\subsection{Wishart and Laguerre processes}

Wishart processes have been introduced in \cite{B89} and \cite{B91}. If $B$ is a 
$n\times n$ Brownian matrix, a Wishart process with parameters $n$ and 
$\delta\geq n+1$
may be obtained as a solution to the matrix-valued SDE
\begin{equation}
 dS_t\,=\, \sqrt{S_t}\,dB_t\,+\,dB_t^{\prime}\sqrt{S_t}\,+\,\delta\,I_n\,dt \;.
\end{equation}
 The eigenvalues process $(\lambda_t^1,\ldots,\lambda_t^n)$ of $\{S_t\}$
satisfies 
\begin{equation}
 d\lambda_t^i\,=\,2\sqrt{\lambda_t^i}\,dW_t^i\,+\,(\delta
 \,+\,\sum_{j\neq i}\frac{\lambda_t^i+\lambda_t^j}{\lambda_t^i-\lambda_t^j})\,dt
 \qquad 1\leq i\leq n \;,
 \end{equation}
and the square roots $r_t^i=\sqrt{\lambda_t^i}$ 
\begin{equation}
 dr_t^i \,=\, dW_t^i \,+\, \frac{1}{2}\frac{\delta-n}{r_t^i}dt\,
   +\frac{1}{2}\sum_{j\neq i}(\frac{1}{r_t^i+r_t^j}\,+\,
  \frac{1}{r_t^i-r_t^j})\,dt
\end{equation}
where $(W^i,\ldots,W^n)$ is a $n$-dimensional Brownian motion. N.Demni \cite{D09}
has remarked that this system is a MSDE with
\begin{equation}
 \Phi(r^1,\ldots,r^n)\,=\, -\frac{1}{2}
[(\delta-n)\sum_i\log \,r^i\,+\,\sum_{i>j}
   \log(r^i+r^j)\,+\,\sum_{i>j}\log(r^i-r^j)]
\end{equation}
on $\{0<r^1<\ldots <r^n\}$ and $\infty$ elsewhere. The system  (53) has a strong 
solution for $\delta>n$. If $\delta=n$, we must add to the right hand 
side of (53) a local time
at $0$ that disappears in (52). It has been proven in \cite{B91} that the eigenvalues 
never collide and if moreover $\delta\geq n+1$ the smallest one never vanishes. 
 This is in accordance with Theorem 6.
 
 Laguerre processes are Hermitian versions of Wishart processes. Only constants are changed in (52), (53) and (54).
 
 \subsection{Reflection groups and Dunkl processes}
 
We only give a short introduction to this topic and refer
to \cite{H90} and \cite{RV98}  for more details. 
For $\alpha\in \R^N\setminus \{0\}$ we denote
by $s_{\alpha}$ the orthogonal reflection with respect to the hyperplane
$H_{\alpha}$ perpendicular to $\alpha$:
\begin{equation}
s_{\alpha}(x)\,=\, x\,-\,2\frac{\alpha.x}{|\alpha|^2}\;.
\end{equation}
A finite subset $R\subset \R^N\setminus \{0\}$ is called a {\sl root system} if
for all $\alpha\in R$
\begin{equation}
   \begin{array}{l}
   R\cap \R \alpha\,=\,\{\alpha,-\alpha\} \;; \\
   s_{\alpha}(R)\,=\,R \;.
   \end{array}
\end{equation}
The group $W\subseteq O(N)$ which is generated by the reflections 
$\{s_{\alpha}, \alpha\in R\}$ is called the {\sl reflection 
group} associated with  $R$. Each hyperplane 
$H_{\beta}:=\{x\in \R^N:\beta.x=0\}$ with 
$\beta\in \R^N\setminus \cup_{\alpha\in R}H_{\alpha}$ 
separates the root system $R$ into $R_+$ and 
$R_-$. Such a set $R_+$ is called a {\sl positive subsystem}
and defines the {\sl positive Weyl chamber} $C$ by
\begin{equation}
 C:= \,\{x\in \R^N:\,\alpha.x>0 \;\; \forall \alpha\in R_+\}\;.
\end{equation}
 A subset $S$ of $R_+$ is called {\sl simple} if $S$ is a vector basis for 
 $span(R)$. The elements of $S$ are called {\sl simple}. Such a subset 
 exists,   is unique and we actually get
 \begin{equation}
  C= \,\{x\in \R^N:\,\alpha.x>0 \;\; \forall \alpha\in S\}\;.
\end{equation}

A function $k:R\rightarrow \R$ on the root system is called a 
{\sl multiplicity function} if it is invariant under the 
natural action of $W$ on $R$. If the multiplicity 
function $k$ is positive on $R_+$, we define the radial 
Dunkl process $X^W$ as the $\overline{C}$-valued continuous path Markov process
whose generator is given by 
\begin{equation}
{\cal L}_k^W u(x)\,=\, \frac{1}{2}\Delta u(x)\,+\, \sum_{\alpha\in R_+}
 k(\alpha)\,\frac{\alpha.\nabla u(x)}{\alpha.x}
 \end{equation}
 for $u\in  C^2(\overline{C})$ with the boundary condition
 $\alpha.\nabla u(x)=0$ for $x\in H_{\alpha}$. Then $X^W$ may be 
viewed as the solution to the MSDE 
\[
  dY_t\,=\,dB_t\,-\,\nabla\Phi(Y_t)\,dt
\]
where $B$ is a Brownian motion and 
\begin{equation}
  \Phi(y)\,=\, \sum_{\alpha\in R_+} k(\alpha)\,\log(\alpha.y)
\end{equation}
on $C$ and $\Phi =\infty$ elsewhere. It was proved in 
(\cite{C06}, \cite{CGY08}) that
this equation has a unique strong solution and if moreover $k(\alpha)\geq 1/2$
for any $\alpha\in R$ then the process never hits the 
walls $H_{\alpha}$ of the Weyl chamber. In \cite{D09}, it is proved that if 
$k(\alpha)<1/2$ for a simple root $\alpha$, then the process 
hits $H_{\alpha}$ a.s. As a consequence of this result and of Theorem 6
(see also the statement at the bottom of p.117 in \cite{CGY08}), 
we are in a position to classify the boundary behavior of the radial
Dunkl process in the Weyl chamber.
\begin{Prop}
For any $\alpha \in R_+$ let 
$\sigma_{\alpha}:=\,\inf\{t>0\,:\,X_t^W \in H_{\alpha}\}$.
\begin{itemize}
\item
 If $\alpha\in R_+ \setminus S$, then $\P(\sigma_{\alpha}=\infty)=1$,
\item
 If $\alpha\in S$ and $k(\alpha)\geq 1/2$, then $\P(\sigma_{\alpha}=\infty)=1$,
\item
 If $\alpha\in S$ and $k(\alpha) < 1/2$, then $\P(\sigma_{\alpha}<\infty)=1$.
\end{itemize}
\end{Prop}   

\subsection{Trigonometric and hyperbolic interactions}

Others interactions have been studied in \cite{CL01}.

The trigonometric system (\cite{D62}, \cite{HW96}, \cite{S98}) reads
\begin{equation}
 \begin{array}{lll}
  dX_t^j\,=\,dB_t^j\,+\,\frac{\gamma}{2}\sum_{k\neq j}
   \cot\frac{X_t^j-X_t^k}{2} & &1\leq j\leq n  \\
    X_t^1 \leq X_t^2 \leq \ldots \leq X_t^n\leq X_t^1\,+\,2 \pi 
 \end{array} 
\end{equation}
This can be interpreted as the solution to the MSDE associated with 
\begin{equation}
  \Phi(x) \,=\, \sum_{i>j}\phi(x.\frac{e_i-e_j}{\sqrt{2}})
  \,+\,\sum_{i<j}\phi(x.\frac{e_i-e_j}{\sqrt{2}}\,+\,\pi\sqrt{2})
\end{equation}
where 
\begin{equation}
 \begin{array}{llll}
  \phi(u) & = & -\gamma\,\log\,(\,\sin\frac{u}{\sqrt{2}}) & 
  0<u<\frac{\pi}{\sqrt{2}} \\
     & = & \infty  & \mbox{elsewhere.}
 \end{array}
\end{equation}    
It has been proved in \cite{CL01} there exist a.s. collisions if $\gamma<1/2$.

The hyperbolic system (\cite{NRW86}, \cite{S07}) is
\begin{equation}
 \begin{array}{lll}
  dX_t^j\,=\,dB_t^j\,+\,\gamma \sum_{k\neq j}
   \coth \,(X_t^j-X_t^k) & &1\leq j\leq n  \\
    X_t^1 \leq X_t^2 \leq \ldots \leq X_t^n \;. 
 \end{array} 
\end{equation}
In this case
\begin{equation}
 \Phi(x)\,=\,\sum_{1\leq j<k\leq n}\phi(x.\frac{e_k-e_j}{\sqrt{2}})
\end{equation}
with
\begin{equation}
 \begin{array}{llll}
  \phi(u) & = & -\gamma\,\log\,(\,\sinh(\sqrt{2} u)) & u>0 \\
     & = & \infty  & \mbox{elsewhere.}
 \end{array}
\end{equation}     
and collisions occur with positive probability if $\gamma<1/2$.

\end{document}